# LAPLACE, FOURIER, AND STOCHASTIC DIFFUSION

T. N. Narasimhan
Dec 13, 2009

*T. N. Narasimhan's interest in mathematics stems from his use of diffusion equation to study water flow, chemical migration, and heat transport in soils, aquifers, and other geological systems. After obtaining a doctorate in engineering science from the University of California at Berkeley, he held joint appointments at Berkeley in the departments of Materials Science and Engineering, and Environmental Science, Policy and Management. He is affiliated with the Earth Sciences Division of Lawrence Berkeley National Laboratory. After three decades of service, he retired in 2006.*

**ABSTRACT**  Stochastic diffusion equation, which attained prominence with Einstein's work on Brownian motion at the beginning of the twentieth century, was first formulated by Laplace a century earlier as part of his work on Central Limit Theorem. Between 1807 and 1811, Fourier's work on heat diffusion, and Laplace's work on probability influenced and inspired each other. This brief period of interaction between these two illustrious figures must be considered remarkable for its profound impact on subsequent developments in mathematical physics, probability theory and pure analysis.

## 1   INTRODUCTION

Stochastic diffusion attained prominence with Albert Einstein's 1905 paper [3] on random motion of colloidal particles in water. Central to Einstein's contribution was an equation, analogous to Fourier's heat equation, with the dependent variable being probability density, rather than temperature. But, the equation had been introduced a century earlier by Pierre Simon Laplace as part of his major work on probability theory, especially what we now call the Central Limit Theorem. During a brief four-year period, 1807-1811, Laplace's path-breaking investigation of probability intersected with Joseph Fourier's profound study of heat movement in solids. This fortunate intersection led to these fertile minds mutually influencing and inspiring each other, resulting in the creation of modern mathematical physics. Although Laplace was a physicist at heart, his pioneering work on stochastic diffusion was of abstract, mathematical nature, while Fourier's work was devoted to bringing a quantitatively observable physical process within the folds of rigorous mathematics. Philosophically, this dichotomy of diffusion, unifying the abstract and the observable is noteworthy, as has been pointed out by Narasimhan [20]. Indeed, Fourier himself was struck by this fascinating consilience. It is fitting to start with Fourier's own words.

In "Preliminary Discourse" of Analytical Theory of Heat, Fourier [7, p. 7] observed[1], *"We see, for example, that the same expression whose abstract properties geometers had considered, and which in this respect belongs to general analysis, represents as well the motion of light in the atmosphere, as it determines the laws of diffusion of heat in solid matter, and enters into all the chief problems of the theory of probability."* Clearly, Fourier was aware that his theory of heat and the theory of probability had mathematical connections. If so, what is the "*qu'une même expression*" that impressed Fourier? What metaphors connect the observable and the abstract?

This paper explores these questions. Fourier's 1807 masterpiece disappeared after his death in 1830, and was discovered more than fifty years later by Gaston Darboux. But, until the publication of a detailed account of this work in 1972 by Grattan-Guinness [8], it went largely unnoticed, and Fourier's 1822 book [6,7] was the principal source of information on his contribution. However, Grattan-Guinness' work brought to light many important developments that occurred between 1807 and 1822 relating to acceptance of Fourier's work by the leadership of French mathematics, including Laplace. Mathematical aspects of these developments have since been addressed by historians Herivel [12], Bru [2], and Gillespie [9]. Additionally, Hald [11] gives a comprehensive account of Laplace's work on probability leading to the publication of his masterpiece *Théorie analytique des probabilites* [17]. The present work complements these contributions by focusing attention on the interactions between Fourier and Laplace between 1807 and the publication of Fourier's book in 1822. The perspectives presented seek a comparative understanding of the beginnings of physical and stochastic diffusion.

In a modern sense, "stochastic diffusion" implies uncertainty associated with random processes. Extension of probability theory to random variables was pioneered during the middle of the nineteenth century by Augustine Cournot who proposed a practical theory of random variables [2]. During a better part of the 18$^{th}$ century, probability was associated with mathematics of actual games, theory of risks, and thought experiments involving hypothetical urns. In a paper presented to the French Academy on March 10, 1773, Laplace made a bold departure, and applied probability to celestial mechanics for studying causes of events. Later, he extended probability to the problem of correcting instrumental error in physical observations [9]. Thus, Laplace's work was devoted to investigation of theory of errors, rather than formal study of random variables. It is worth noting that the phrase Central Limit Theorem came to be established only during the twentieth century, at the suggestion of Polyà [21, 23].

---

[1]On voit, par exemple, qu'une même expression, dont les géomètres avaient considéré les propriétés abstraites et qui, sous ce rapport, appartient à l'Analyse générale, représénte aussi le mouvement de la lumière dans l'atmosphère, qu'elle détermine les lois de la diffusion de la chaleur dans la matière solide, et qu'elle entre dans toutes les questions principales de la Théorie des probabilités [Fourier, 6, p. xxiii]

## 2 CLUES TO CONNECTIONS BETWEEN DIFFUSION AND PROBABILITY

Spreading of heat in a solid is governed by combined effects of thermal conductivity and thermal capacity of a solid material, as heat is driven from a location of higher to one of lower temperature. Thermal conductivity is a proportionality constant linking heat flux with temperature gradient, and thermal capacity, first defined and measured by Lavoisier and Laplace [19], quantifies the relationship between magnitude of temperature and quantity of heat stored in a body. Spreading of probability represents increase in uncertainty in cumulative error of a finite number of random variables in proportion to number of variables summed up.

Fourier's [4] parabolic equation for heat diffusion is,

$$K \frac{\partial^2 T}{\partial x^2} = C \frac{\partial T}{\partial t},$$  (1)

where K is thermal conductivity, T, temperature, C, thermal capacity, x, distance along the abscissa, and t, time. A diffusion problem is fully defined when (1) is augmented by appropriate boundary and initial conditions. In his 1807 monograph, Fourier devoted attention exclusively to bounded, symmetrical solids [rod, prism, sphere, cube, ring]. Assuming K and C to be independent of temperature, and thus linearizing the differential equation, he pioneered a number of novel mathematical techniques for solving it.

In stochastic diffusion, spreading is quantified by increase in variance of cumulative error distribution in direct proportion to number of samples. The Central Limit Theorem is a mathematical elaboration of this fact. Let $x_1, x_2, \ldots x_n$ be independently and identically distributed random variables with frequency function f(x), mean µ, and variance $\sigma^2$. The Central Limit Theorem states that $s_n = x_1 + x_2 + \ldots + x_n$ asymptotically approaches normal distribution with mean nµ and variance $n\sigma^2$ [11].

Comparison of the two phenomena shows that time in the thermal process is analogous to number of samples in error propagation. However, the error propagation problem has no feature analogous to the bounding surface of a solid in the heat flow problem. Consequently, to understand how the "same expression which determines the laws of heat diffusion in solid matter enters into all the chief problems of the theory of probability", it is necessary to examine how Fourier posed and solved heat flow problems in infinite media. Looking carefully at Fourier's Analytical Theory of Heat [6,7] from this perspective, the connection is readily found in Article 364, which deals with transient heat flow along an infinite line, over a segment of which arbitrary initial conditions are prescribed, with zero temperature elsewhere.

## 3 LAPLACE (1809), AND FOURIER (1811)

Laplace, who had actively worked on probability theory for over fifteen years from 1771, devoted most of his attention over the next twenty years to the study of planetary mechanics, culminating in the publication of *Traité de mecaniqué céleste* between 1799 and 1805. He returned to probability thereafter, and provided the first proof of the Central Limit Theorem in 1810 for variables with a continuous uniform distribution [11]. Leading up to this work, in 1809, Laplace [14] investigated the use of génératrice functions to evaluate the probability that the sum of a given number of identically distributed random variables would take on a given value. It was well known through earlier investigations of James Bernoulli, Abraham de Moivre and others that the required probability constituted the coefficient of a particular term in a power series stemming from the génératrice function. The difficulty was that estimating the numerical magnitude of the coefficient was mathematically very difficult when the number of random variables summed up became very large. Joseph Lagrange and Laplace devoted much of their energies to evaluating the required coefficients using finite difference equations and associated definite integrals.
In the subsection entitled, "On definite integrals of partial difference equations"[2], Laplace [14, p. 235] started with the study of linear second-order partial differential equations whose solution involved two arbitrary functions. In a special case, this general problem reduced to one whose solution involved but a single arbitrary function. This was a parabolic equation of the form, in Laplace's notation [14, p. 238],

$$\frac{ddu}{ds^2} = \frac{du}{dx'}. \qquad (2)$$

---

[2] Sur les Intégrales défines des Équations à différences partielles

He then showed how the estimation of the coefficient of a particular term in a power series expansion could be transformed to finding a solution to the aforesaid equation.

Let u be a power series expansion in t and t', with $y_{x,x'}$ being the coefficient of ($t^{x.}\, t'^{x'}$.) Let u be the génératrice function for $y_{x,x'}$, and $u[[1/t - 1]^2 - [1/t' - 1]]$ be the génératrice function for $\Delta^2 y_{x,x'} - \Delta' y_{x,x'}$, the characteristic $\Delta$ being relative to x and $\Delta'$ being relative to x'. The goal was to solve for the coefficient $y_{x,x'}$, when the number of terms to be considered in the power series are large. To this end, Laplace used recursive relations to arrive at the finite difference equation,

$$\Delta^2 y_{x,x'} = \Delta' y_{x,x'}. \qquad (3)$$

In the infinitesimal limit, this led to the differential equation [in Laplace's notation],

$$\frac{d^2 y}{dx^2} = \frac{dy}{dx'}. \qquad (4)$$

Based on his earlier work, Laplace then demonstrated that,

$$y = \int_{-\infty}^{+\infty} c^{-z^2}\, \varphi[x + 2z\sqrt{x'}]\, dz. \qquad (5)$$

satisfied the partial differential equation (4), where $\varphi$ is any arbitrary function. Here, $y_{x,x'}$ represents the probability that the sum of x' identically distributed random variables takes on the value x. Comparing with the heat equation, probability y corresponds to temperature, the magnitude of the sum of random variables, x, corresponds to distance x, and the number of random variables, x', corresponds to time. The coefficient, $y_{x,0}$, represent initial conditions.

The form of Laplace's solution immediately revealed to Fourier that he could seek integral solutions to the heat equation in addition to the series solutions. Thus inspired, Fourier [5] expanded his 1807 work, and filed it with the Institut de France on September 28, 1811 in response to the prize competition it had set up. The significant addition in this expansion was Chapter XI on the linear movement and variation of heat in a body with one infinite dimension. For the first time, Fourier addressed the problem of heat movement in a solid without a bounding surface. Problems of this type are driven solely by initial conditions.

In particular, Fourier considered an infinite line with $-\infty < x < +\infty$. At time t = 0, the temperature everywhere along this line was zero, except over a segment extending on either side of x = 0. Over the segment, temperature distribution was an arbitrarily prescribed function f(x). He considered several cases with f(x) representing different patterns of temperature variation over the segment. The governing differential equation was,

$$\frac{du}{dt} = k \frac{d^2u}{dx^2}, \qquad (6)$$

with initial condition f(x). To solve for u(x,t), he sought solutions in three different forms, two of them involving convolution integrals with the heat kernel,

$$u = \int e^{-kq^2 t} f(x) \cos qx \, dq, \qquad (7a)$$

$$u = \int e^{-kq^2 t} f(x) \sin qx \, dq, \quad \text{and} \qquad (7b)$$

$$u = e^{-x} e^{-kt}. \tag{7c}$$

Fourier went on to show through a series of transformations that (7c) yielded solution of the form,

$$u = \int_{-\infty}^{+\infty} c^{-z^2} \varphi[x + 2z\sqrt{x'}] \, dz. \tag{8}$$

Note that , (8) has the same form as Laplace's solution (5) to the probability problem (4).

Although Fourier [5] did not refer to Laplace's work in the Prize Essay, he acknowledged in Fourier [7, Art. 364][3], *"This integral which contains one arbitrary function was not known when we had undertaken our researches on the theory of heat, which were transmitted to the Institute of France in the month of December, 1807: it has been given by M. Laplace, in a work which forms part of Volume VIII of the Mémoires de l'École Polytechnique; we apply it simply to the determination of the linear movement of heat"*. With some modifications and change of symbols these results were presented in Fourier [6,7] as Chapter IX, Section I. His conclusion at the end of this section was that solutions to equation (6) arrived at through different forms [e.g. (7a), (7b), (8)] were equivalent.

Following his 1809 contribution, Laplace chose to pursue proof of the Central Limit Theorem using characteristic functions rather than the differential equation, and announced his result to the Academy in April 1810 [2,22]. Soon thereafter, he found that the appearance of Carl Friedrich Gauss' recently published work on the method of least squares had clearly shown the connection between Central Limit Theorem and linear estimation [22].

The following year he generalized a two-urn problem of Bernoulli and formulated a model involving what is now referred to as Markov Chain with transition probabilities. This model led to a second order partial differential equation,

---

[3] Cette intégrale, qui contient une fonction arbitraire, n'était point connue lorsque nous avons entrepris nos recherches sur la Théorie de la chaleur, qui ont été remises à l'Institut de France dans le mois de décembre 1807; elle a été donnée par M. Laplace , dans un Ouvrage qui faite partie du Tome VIII du *Journal de l'École Polytechnique* [[1]]; nous ne faisons que l'appliquer à la détermination du mouvement linéaire de la chaleur. [Fourier, 6, p. 414]

$$\frac{\partial u}{\partial t} = 2u + 2\mu\frac{\partial u}{\partial \mu} + \frac{\partial^2 u}{\partial \mu^2}. \tag{9}$$

His solution to this problem involved polynomials, which would subsequently be recognized as being proportional to Hermite polynomials [11]. The solution anticipated Fourier-Hermite series for functions defined over infinite domains.

## 4  THE EXPRESSION $\dfrac{e^{\frac{-x^2}{4t}}}{\sqrt{t}}$

In Analytical Theory of Heat, Fourier [6,7] considered heat movement in infinite solids. He started with $u = e^{-n^2 t} \cos nx$ as satisfying the differential equation, $\dfrac{du}{dt} = \dfrac{d^2 u}{dx^2}$, and showed that

$$\int_{-\infty}^{\infty} e^{-n^2 t} \cos nx \, dn = e^{\frac{-x^2}{4t}}. \tag{10}$$

Consequently, the aforesaid differential equation is satisfied by,

$$u = \frac{1}{2\sqrt{\pi}} \int_{-\infty}^{\infty} f(\alpha) \frac{e^{\frac{-(x-\alpha)^2}{4t}}}{\sqrt{t}} d\alpha, \tag{11}$$

where α is any constant. If we let $(x - \alpha)^2 / 4t = q^2$, then,

$$u = \frac{1}{\sqrt{\pi}} \int_{-\infty}^{\infty} e^{-q^2} f(x + 2q\sqrt{t})\, dq .\tag{12}$$

If, in (11), we set t = σ², then, $\dfrac{1}{\sqrt{2\pi\sigma^2}} e^{\frac{-(x-a)^2}{4\sigma^2}}$ is the probability density function for norma distribution with mean α and variance σ². Thus, Fourier established that probability density function which plays a fundamental role in probability theory, also forms part of solutions fundamental to transient heat diffusion in infinite media.

## 5  MOVEMENT OF LIGHT IN THE ATMOSPHERE

Now we consider Fourier's [7, p. 7] reference to, "...*the same expression .... represents as well the motion of light in the atmosphere, ...*"[4]. Presumably, Fourier was referring to Laplace's work, "*Mémoire sur les mouvements de la lumière dans les milieux diaphanes*", read before the Académie in 1808, and published in 1810 [16].

A major part of this work was devoted by Laplace to developing a theory for transmission of light in transparent media, including the atmosphere, based on the philosophy of action-at-distance. Laplace believed so strongly in this philosophy [10, p. 93] that he extended it, by analogy with light propagation in the atmosphere, to a discussion of the movement of heat in solids. This he did in a long "Note". The central concept in this approach was that of interacting molecules of light, or analogously, of heat. Laplace prefaced his discussion of heat propagation with the statement, "*By considering action-at-distance of molecule and molecule, and extending such action to heat, we arrive, through a simple and precise way, at the true differential equations that*

---

[4] "…..qu'une même expression, ….. représente aussi le mouvement de la lumière dans l'atmosphère, ..." [Fourier, 6, p. xxiii]

*describe heat movement in solid bodies and its variations on their surface, and thus this very important branch of physics enters in the area of Analysis"*.[5]

For Laplace, action-at-distance in regard to heat was embodied in Newton's principle that the quantity of heat communicated by a body to its neighbor is proportional the difference in their temperatures. Based on this, he then presented the partial differential equation for the flow of heat in a solid by analogy with similar derivation for the propagation of light [16, p. 293; Laplace's notation],

$$du = a \, dt \, \frac{\partial^2 u}{\partial x^2}, \tag{13}$$

where the constant a is thermal conductivity. He then observed that this equation can be generalized to three spatial dimensions.

It is not clear if this Note was prepared in 1808 when the paper was presented before the Institut, or it was prepared in 1810. Regardless, it is clear that Laplace implicitly conceded Fourier's priority in presenting the parabolic equation. However, his desire seems to be one of providing a better way of deriving the equation than what had been achieved by Fourier. This is evident in his assertion, "*However, just like mathematicians arrived at the equations describing the movement of light in atmosphere starting from an inaccurate hypothesis, the hypothesis that the action of heat is limited to contact area can lead to the equations describing heat movement inside and at the surface of bodies. I need to take note that M. Fourier already arrived at these equations, the real bases of which seem to be those I just presented.*"[6]

---

[5] Enfin la considération des actions *ad distans* de molécule à molécule, étendue à la chaleur, conduit d'une manière claire et précise aux véritables équations différentielles du mouvement de la chaleur dans les corps solides et de ses variations à leur surface, et par là cette branche tres important de la Physique rentre dans le domaine de l'Analyse [16, p. 290].

[6] Mais, de meme que les geometres avaient été conduits aux equations du mouvement de la lumière dans l'atmosphère, en partant d'une supposition inexacte, de même l'hypothèse de l'action de la chaleur limitée au contact peut conduire aux équations du mouvement de la chaleur dans l'intérieur et à la surface des corps. Je dois observer

Page 10 of 17

From the foregoing, there is little doubt that Fourier's statement, "*.... represénte aussi le mouvement de la lumière dans l'atmosphère.....*" refers to Laplace [16]. Given that, it is pertinent to examine how Fourier approached the derivation of the same heat equation. As has been described in detail by Grattan-Guinness [8], and Herivel [12], Fourier began his investigation of heat around 1804, starting with action-at-distance and Newton's principle. In this, he followed the same line of reasoning as Biot [1] before him. However, he encountered difficulties in formally setting up a differential equation. Consequently, he abandoned action-at-distance, and introduced the continuity assumption that the state of heat at a point depends solely on the immediately preceding point. This approach essentially introduced the notion of a continuum. That this approach has withstood the test of time suggests that Laplace's claim that his method provides the "real base" of the heat equation does not carry conviction.

## 6     PERSONAL INTERACTIONS

As we have seen, the period 1807 - 1811 was remarkable in the history of mathematical statistics and mathematical physics. On the human side, this period was distinguished by an initial rivalry and subsequent rapprochement between two intense individuals who were revolutionizing science. What was the nature of this personal interaction?

There is little doubt that Fourier was the first to formulate the parabolic equation in 1807. Although he had experimented for three decades with a variety of transforms (including what would later be termed as Fourier transform) to solve difference equations and differential equations, Laplace apparently did not recognize that the parabolic equation would help evaluate not only the mean but also the variance as well of the sum of a large number of random variables [22]. His 1809 formulation was clearly catalyzed by Fourier's 1807 monograph. On his part, Fourier was inspired by Laplace's 1809 formulation of the parabolic equation to recognize that heat flow in infinite domains constituted a new class of problems, and that solutions to the heat equation can also be obtained in the form of integrals. Clearly, physical diffusion and stochastic diffusion had mutually influenced each other at birth.

---

que "M. Fourier est déjà parvenu à ces équations, dont les véritables fondements me paraissent être ceux que je viens de présenter. [16, p. 295]

During the period 1807 to 1811, when both Laplace and Fourier were intensely addressing their respective topics, there was some tension between them. While Laplace [14], in his 1809 work, did not acknowledge Fourier's 1807 work, Fourier [5] failed to cite Laplace's work in his Prize Essay. Laplace [16, p. 295] conceded that Fourier had already presented the heat equation, but asserted that his derivation based on action-at-distance was more fundamental than Fourier's derivation. The tension gradually gave way to mutual respect when Fourier spent nearly a year in Paris, starting from the summer of 1809. During this stay, Fourier regularly attended meetings in Laplace's estate at Arcueil, the uncontested center of world science at that time [2,8]. Thus, although Fourier [6, p. xxiii] did not specifically mention Laplace, it is clear that he was referring to Laplace in stating, "...*qu'une mème expression, don't les géométres avaient considéré les propriétés et qui, sous ce rapport, appartient à l'Analyse générale,...*".  For his part, Laplace [18, p. 83] reproduced the heat equation and the equation at the boundary and complimented Fourier by stating, "... *M. Fourier was the first to present the fundamental equations (1) and (2) in the excellent paper that won the prize proposed by the Institute on the Theory of Heat; I shall give their demonstration in a different book.*"[7] Presumably, *"j'en donnerai la démonstration dans un autre livre"* refers to Laplace [16, p. 295].

As cited at the beginning of this paper, Fourier was impressed by the fact that the "same expression" which determines the laws of diffusion of heat also enters into all "chief problems of the theory of probability". What "same expression" was Fourier fascinated about? From what we have seen, there are three possibilities. The first is the convolution integral,

$$\int e^{-q^2} f[x + 2q\sqrt{t}]\, dq ,$$

which he borrowed from Laplace's 1809 paper. The second expression,

$$\frac{e^{\frac{-(x-a)^2}{4t}}}{\sqrt{t}} ,$$

which is essentially the same as the probability density for normal distribution. The third is the parabolic equation itself, in view of Laplace's work on propagation of light in the atmosphere.

---

[7] M. Fourier a donné le premier les équations fondamentales (1) and (2) dans l'excellente pièce qui a remporté le prix proposé par l'Institut sur la Théorie de la chaleur; j'en donnerai la démonstration dans un autre livre [Laplace, 1823, p. ??]

## 7      CONCEPT OF A FUNCTION AND PHYSICAL IMPLICATIONS

To Lagrange, Fourier's statement that an arbitrary function could be expressed as a trigonometric series was so unexpected that he opposed it strongly.  In a recent paper, Kahane [13] presents new evidence on Lagrange's erroneous criticism, based on a "Schriftstück" of Lagrange mentioned in Bernhard Riemann's Habilitation dissertation. As Kahane shows, Lagrange's criticism was an indication that the concept of a function, which was a source of controversy among Jean d'Alembert, Leonhard Euler, and Daniel Bernoulli during the $18^{th}$ century, was still evolving around 1800.  Indeed, Fourier's work inspired Augustine Cauchy, Lejeune Dirichlet and Riemann to continue to refine the concept of a function to pave the way for Georg Cantor and others to lay the foundations of modern theory of functions of a real variable.

Lagrange's criticism stemmed from the fact that he took a trigonometric series defined over $0 < x < \pi/2$ and showed that it led to an inconsistent result when x was set to 0.  In response, Fourier had pointed out that equations of certain type cannot be used without specifying the limits between which the values of the variable have to be considered.  But Lagrange did not relent. Fourier's contribution to the theory of functions was to establish that a function is only valid over a specified domain.

Against this background it is worthwhile to examine the physical implications of functions and their domains. Whereas Fourier's 1807 monograph established that any function defined over a bounded domain could be represented by trigonometric series, Laplace's solution of a general second-order partial differential equation (9) in terms polynomials anticipated the later development showing that an arbitrary function over an infinite domain could be represented using Fourier-Hermite series. Thus, finite and infinite domains fall into distinct categories in terms of representative functions.

Finite and infinite domains also relate to distinct categories of physical problems.  In a transient system, heat flow is driven by non-uniform spatial distribution of temperature at the initial time, or by external influences acting on the bounding surfaces of the system, or both.  In infinite, unbounded solid bodies, the initial condition is the sole cause of heat flow.  The self-smoothing tendency of the system is to dissipate disturbances by itself, without any external influence.  Here, the fundamental problem of interest in an infinite system is the release of a certain amount of heat in the vicinity of a point at time zero, and to predict the spreading (diffusion) of heat for t > 0. This is referred to as an instantaneous source. All other problems pertaining to an infinite system can be solved by superposition of this fundamental solution using convolution integrals. The fundamental problem is inherently symmetrical. Laplace's stochastic diffusion problem is also an initial value problem. Given a probability density distribution for n = 0, n being analogous to time, one solves for the spreading of probability density as n becomes progressively large. In the process, no restrictions are placed on the randomness of the values sampled. Time in the physical problem and number of samples in the stochastic problem are unbounded and tend to infinity. Therefore, in problems involving infinite domains, the time derivative, or equivalently, the derivative with reference to number of samples n has to be non-zero. An equilibrium state is not theoretically definable.

It is interesting that recursive formulas and difference equations played a very important role in many of probability problems solved by Laplace. Laplace has remarked [11, p. 338] that if the initial distribution were known, all subsequent distributions can be calculated with the help of the recursive formula. This remark reinforces the view that evolution in time (or, equivalently, number of samples) is central to stochastic diffusion.

In contrast, boundary-value problems are driven by forces imposed on the boundary by external causes. Under time-invariant boundary conditions, the system is driven to steady state flow characterized by vanishing time derivative. The solution satisfies Dirichlet Principle, an integral that has to be minimized. The unique solution is independent of any initial condition that may have existed at the beginning.

In stochastic diffusion, spreading continues with progressively increasing number of samples as long as random sampling continues unfettered by external influence. Therefore, if any bounds on are set on the value of the sum of the random variables, then randomness is inhibited by external causes. If sampling is continued under bounded conditions, system progress will be influenced more and more by boundary conditions.

Thus, problems involving infinite domains and finite domains constitute two distinct classes. Mathematically, the behavior of the latter can be described using trigonometric series, while the former can be described using Hermite polynomials. Initial-boundary value problems may be considered to be mixed problems, combining features of both. Semantically, it is interesting to note that in Analytic Theory of Heat, Fourier [6,7] uses the word "diffusion" in the heading of Chapter IX devoted to infinite media. This is eminently reasonable because "diffusion" or "spreading" can occur only in infinite domains where no external forces inhibit spreading. Fourier himself noticed this difference between problems defined over finite and infinite domains when he [7, Article 343] stated[8], "*In the problems we previously discussed, the integral is subjected to a third condition which depends on the state of the surface: for which reason the analysis is more complex, and the solution requires the employment of exponential terms. The form of the integral is very much more simple, when it need only satisfy the initial state;...*"

## 8  CONCLUDING REMARK

Laplace and Fourier were natural philosophers seeking to comprehend a finite world subject to errors of discrete observations. Their difference equations and recursive relations were useless when the number of observations were large. To overcome this difficulty, their creative intellects led them from difference equations to differential equations and a host of definite integrals and convergent algebraic series. Yet, the observational world remains finite and discrete, and the algebraic expressions are but idealized approximations of reality. We continue to grapple with balancing the discrete and the continuous. An intriguing question emerges: if a digital computer had been available to Lagrange, Laplace, and Fourier to handle large numbers, what course would mathematics have taken?

## 9  ACKNOWLEDGMENTS

I am very grateful to Jean-Pierre Kahane for many illuminating discussions on history of mathematics in France as well as for deep insights into the theory of functions and its relation to mathematical physics. I thank Ghislain de Marsily, and Roger Hahn for criticisms and valuable suggestions. Norbert Schappacher's critical review helped greatly in revising the paper. Over the past decade, I have had many valuable exchanges with Stephen Stigler on the history of mathematical statistics.  These exchanges have been of help in preparing this paper.

---

[8] Dans les questions que nous avons traitées précédemment, l'integrale est assujettie à une troisième condition qui dépend de l'état de la surface. C'est pour cette raison que l'analyse en est plus composée et que la solution exige l'emploi des termes exponentiels. La forme de l'integrale est beaucoup plus simple lorsqu'elle doit seulement satisfaire à l'etat initial, ... . [Fourier, 6, p. 388].

**REFERENCES**

[1]     Biot, J. B.;  *Mémoire sur la chaleur*, *Bibliotheque Britannique*, 27, 310-329, 1804.
[2]     Bru, B.; Poisson, the probability calculus, and public education. English translation of an article published in 1981 by Glenn Shafer, *Journ@l Electroniqued'Histoire des Probablités et de la Statistique*, 1 (2), 25 p.,  November 2005.
[3]     Einstein, A.;  Über die von der molekularkinetischen Flussigkeiten suspendierten Teilchen. *Ann. Der Phys.*, 17 (1905), 549-560, 1905.  [English translation in *Investigations on the theory of Brownian movement by Albert Einstein,* Edited with notes by R. Fürth, pp. 1-35, Methuen, London, 1926]
[4]     Fourier, J. B. J.;  *Mémoire sur la propagation de la chaleur*. Read before the Institut on December 20, 1807.  Unpublished, and preserved in MS 1851 of the École Nationale des Ponts et Chaussées, Paris.  Full text with commentary by I. Grattan-Guinness [8].
[5]     Fourier, J. B. J.; Théorie du mouvement de la chaleur dans les corps solides, The Prize Essay, deposited with the Institut on September 28, 1811. Published in 1824 as *Mémoires l'Académie Royale des Sciences de l'Institut de France, years 1819 and 1820*, 185-556, with one Plate containing 13 figures.
[6]     Fourier, J. B. J.; *Théorie Analytique de la Chaleur*. C. F. Didot, Paris, 563 p., 1822a.
[7]     Fourier, J. B. J.; *Analytical Theory of Heat*, 1822b. Translated with notes by A. Freeman, the Cambridge University Press, London, 1878, 466 pages.
[8]     Grattan-Guinness,  I., in collaboration with J. R. Ravetz, ;  *Joseph Fourier, 1768-1830*, The M.I.T. Press, Cambridge, Mass., 516 p., 1972.
[9]      Gillespie, C. C., *Pierre Simon Laplace 1749-1827, A life in Exact Science*, Princeton University Press, Princeton New Jersey, 322 p.,1997.
[10]    Hahn, R.; *Pierre Simon Laplace 1749-1827, A Determined Scientist*, Harvard University Press, Cambridge, Mass., 310 p., 2005.
[11]    Hald, Anders.;  *A History of Mathematical Statistics from 1750 to 1930*. John Wiley and Sons, New York, 793 p., 1998.
[12]    Herivel, J.; *Joseph Fourier, The Man and the Physicist*, Clarendon Press, Oxford, 1975.

[13] Kahane, J-P., Partial differential equations, trigonometric series, and the concept of a function around 1800: a brief story about Lagrange and Fourier, *Proc. Symposia in Pure Math.*, 79, 187-205, 2008.

[14] Laplace, P. S.; Mémoire sur divers points d'Analyse, *Journal de l'Écolé Polytechnique*, Tome VIII., 229-265, 1809, Oeuvres Completes, XIV, 178-214.

[15] Laplace, P. S.; Mémoires sur les approximations des formules qui sont fonctions de très grands nombres et sur leur application aux probabilités, *Mémoires l'Académie Royale des Sciences*, Paris, 353-415, 1810, Oeuvres Completes, XII, 301-353.

[16] Laplace, P. S.; Mémoire sur les mouvements de la lumière dans les milieux diaphanes, *Mémoires de l'Académie des Sciences,* I$^{st}$ Série, Tome X, 1810, *Oeuvres Completes de Laplace*, Vol 12, 267-298.

[17] Laplace, P. S.; *Théorie analytique des probabilites*. V Courcier, Paris, 1812, 1814, 506 p.

[18] Laplace, P.S.; De la chaleur de la Terre et de ladiminution se la durée du jour par son refroidissement, 1823, *Oeuvres Completes de Laplace*, Vol 5, p. 83.

[19] Lavoisier, A. L. and Laplace, P. S.; Mémoire sur la Chaleur, *Mémoires de l'Académie Royale des Sciences*, Paris, for the year 1780, 355-408. Paper read June 28, 1783.

[20] Narasimhan, T. N.; On the dichotomous history of diffusion, *Physics Today*, 62(7), 48-53, 2009.

[21] Polyà, G.; Uber den zentralen Grenzwertsatz der Wahrscheinlihkeitsrechung und das Momentenproblem , *Math. Zeit.*, 8, 171-181, 1920.

[22] Stigler, S. M.; *The History of Statistics*, Harvard University Press, 410 p., 1986.

[23] Stigler, S. M.; Personal communication, 2008.

T. N. Narasimhan

Department of Materials Science and Engineering

210 Hearst Memorial Mining Building

University of California

Berkeley, Ca 94720-1760

Email: tnnarasimhan@LBL.gov